\documentclass[12pt,reqno]{amsart}
\usepackage{graphicx}
\usepackage{amsmath}
\usepackage{amssymb}
\usepackage{amscd}
\usepackage{mathrsfs}
\usepackage[all,cmtip]{xy}
 \usepackage{color}
\usepackage{stmaryrd}

\newtheorem{theoremc}{Theorem}

\newtheorem{rk}[theoremc]{Remark}

\newcommand\com[1]{}
\newcommand\C{{\mathbb C}}
\newcommand\Cc{{\let\mathcal\mathscr\mathcal C}}

\newcommand\E{\mathcal{E}}

\newcommand\g{{\frak g}}

\newcommand\op[1]{\mathop{\rm #1}\nolimits}

\newcommand\p{\partial}

\newcommand\R{{\mathbb R}}

\begin{document}

 \title{Symmetries of second order ODEs}
 \author{Boris Kruglikov}
\address{Department of Mathematics and Statistics, UiT the Arctic University of Norway, Troms\o\ 90-37, Norway.
\ E-mail: {\tt boris.kruglikov@uit.no}. }
 \date{}

 \vspace{-15pt}
 \begin{abstract}
It is demonstrated that point symmetry algebras of general analytic second order ODEs, not 
necessary of principal type, can have all dimensions between 0 and 8 except for 7. For the symmetry 
dimension 8 the ODE must be locally trivializable.
 \end{abstract}



 \maketitle

 \vspace{-18pt}

\section{Introduction and the main result}\label{S1}

Second order scalar ordinary differential equations $y''=f(x,y,y')$ have been studied from 
the geometric viewpoint by Sophus Lie. 
He formulated the problem of their local classification with respect to the diffeomorphism group
of the space $\R^2(x,y)$ of independent and dependent variables (point transformations) that was solved
by his student Arthur Tresse \cite{T}, see also \cite{K}. 
In this work Tresse also shows that the symmetry algebra of such an ODE
can have dimensions only 8, 3, 2, 1 or 0; see also \cite{KT} for an independent and more rigorous approach.

More general equations $F(x,y,y',y'')=0$ not resolved with respect to the second derivative arise in
the study of path geometries with singularities and are important in applications to special functions.
Indeed, Bessel, Legendre and hypergeometric differential equation due to Euler and Gauss are special instances
of this type. All Painlev\'e transcendents are also particular cases. To our knowledge the symmetry analysis of this 
general class has never been performed, and the main goal of the present note is to address it.

In what follows we assume $F$ to be analytic and we consider the algebra $\g$ of 
analytic point symmetries of this implicit ODE. The statements below apply to the case, 
when $\g$ consists of vector fields defined in a fixed 
connected neighborhood $U\subset\R^2(x,y)$, also with a specified connected neighborhood $\hat{U}\subset\E=\{F=0\}$ over it,
as well as to the case for germs of the fields from $\g$ at a fixed point $(x_0,y_0)\in U$.

 \begin{theoremc}
Let $F(x,y,y',y'')=0$ be an analytic ODE of the second order, meaning that $F_{y''}\not\equiv0$. 
Then the point symmetry algebra $\g$ can only have dimensions $8,6,5,4,3,2,1,0$ and all of those are realizable.

Moreover $\dim\g>6$ implies that the equation is locally trivializable, i.e. equivalent to $y''=0$ in
a neighborhood of any point.
 \end{theoremc}

Even if the ODE domain $U\subset\R^2$ is connected and simply-connected, a global trivialization
may not exist, as the sraightening diffeomorphism $U\to\R^2$ can only be an immersion but not an embedding.

 \begin{rk}
If $F_{y''}\equiv0$ we get an ODE of the first order. As a rule such equations have an infinite-dimensional 
symmetry algebra. Indeed, this is always the case with ODEs of the principal type $y'=h(x,y)$;
for instance the trivial equation $y'=0$ has the infinite-dimensional symmetry algebra consisting of the vector fields
$\{a(x,y)\p_x+b(y)\p_y\}\subset\mathcal{D}(\R^2)$. However equations not resolved with 
respect to the derivative can have smaller symmetry algebras; for example, the ODE $y=(y')^m$ 
has two symmetries $\p_x,\frac{m-1}m\,x\,\p_x+y\,\p_y$ for $m\ge3$.
 \end{rk}

Let us remind that second order ODEs of the principal type are instances of parabolic geometries of type 
$(SL_3,P_{1,2})$, see \cite{CS}. A complex relative of this is the Levi-nondegenerate CR-geometry in 
dimension 3, which has parabolic type $(SU_{1,2},P_{1,2})$. Indeed, these two geometries have the same complexifications, and 
the passage between them is as follows: for a real hypersurface $\varphi(z,\bar{z})=0$ in $\C^2$
the Segre family $Q_\zeta=\{z\in\C^2:\varphi(z,\bar{\zeta})=0\}$ defines a 2-parametric family of solutions to a
second order ODE (obtained by elimination of $\zeta$).
In \cite{C} E.\,Cartan used his equivalence method to establish, in particular, that symmetry algebras of 
CR-structures in 3D have dimensions 8, 3, 2, 1 or 0, which corresponds to the Tresse result mentioned above.

CR-structures in 3D without Levi-nondegeneracy constraint were considered in \cite{KS}, and it was proved
that if CR-manifold $M^3$ is not Levi-flat then its symmetry algebra has dimension 8, 5, 4, 3, 2, 1 or 0.
This was done by exploiting the above Segre correspondence and studying the respective Fuchsian type second 
order equations. A much simpler proof of what Kossovskiy and Shafikov called the Poincar\'e dimension 
conjecture was presented in \cite{IK}. It used only a combination of some fundamental facts of Lie theory and CR-geometry. 

In this paper we prove our main result in the same spirit. Note however that second order ODEs 
and CR-structures in 3D are analogous when non-degenerate, but not in general. In fact, 
CR-geometry is given by a complex structure 
on a 2-distribution that can fail to be contact in general, but degenerations of path geometries in 3D
consist of collapsing two line fields, spanning a contact distribution
at regular points, to one line at non-resolvable points of implicitly given ODEs, and also one or both of these
lines can blow up to a plane at singular points.

\section{The proof}\label{S2}

Geometrically, the ODE is encoded as 3-dimensional submanifold $\E=\{F(x,y,y_1,y_2)=0\}$ in the space of
2-jets $J^2(\R)=\R^4(x,y,y_1,y_2)$. The latter is equipped  with the Cartan distribution 
$\mathcal{C}_2=\op{Ann}(dy-y_1dx,dy_1-y_2dx)$. This induces line field with singularities
$\mathcal{C}_\E=\mathcal{C}_2\cap T\E$ on $\E$. Another line field with singularities $V$ is tangent 
to the fiber of the projection $\E\to J^0(\R)=\R^2(x,y)$. 
At singular points the line $\langle\p_u\rangle$ is contained in both $\mathcal{C}_\E$ and $V$, 
and each of those can become a 2-plane. These two fields span the distribution $\Pi=\mathcal{C}_\E+V$
on $\E$. At the regular points, where the map $\pi_{2,1}:\E\to J^1(\R)$ is a submersion, i.e.\ 
$\E=\{y_2=f(x,y,y_1)\}$, this distribution pulls back the contact distribution 
$\mathcal{C}_1=\op{Ann}(dy-y_1dx)$ on $J^1(\R)=\R^3(x,y,y_1)$. 
Point symmetries of the ODE are bijective with symmetries of the pair of line fields with singularities
$(\mathcal{C}_\E,V)$ on the equation-manifold $\E$.

Let us note that due to the assumption $F_{y_2}\not\equiv0$ the set of points, where the equation 
has principal type is open and dense. Near each of them $(\E,\mathcal{C}_\E,V)$ has at most 8-dimensional symmetry algebra, and so $\g$
satisfies $\dim\g\leq8$. Moreover by the result of Tresse mentioned above \cite{T,K,KT} the inequality $\dim\g>3$
implies $\g\subset\mathfrak{sl}_3$.

We begin with the proof of the second part of the theorem, namely that if dimension of the symmetry algebra 
is at least 7, then the ODE is locally trivializable. Since $\mathfrak{sl}_3$ has no subalgebras of dimension 7,
we have $\g=\mathfrak{sl}_3$ in this case. Thus for regular points $a\in\E$, where our ODE is of principal type, 
it is flat, i.e. trivializable. In these points the isotropy algebra $\g_a=\mathfrak{p}_{1,2}$ is the minimal
parabolic (Borel) subalgebra in $\mathfrak{sl}_3$, and a neighborhood of $a$ is homogeneous of the type $SL_3/P_{1,2}$.

Clearly the isotropy has dimension $\dim\g_a\in[5,8]$ and the case $5$ yields a locally homogeneous space.
Consider the other cases. If $\dim\g_a=8$, then $a$ is a fixed point of the simple Lie algebra $\g=\mathfrak{sl}_3$
and so, by the Guillemin-Sternberg theorem \cite{GS}, the algebra linearizes at $a$. 
However for the standard representation
of $\mathfrak{sl}_3$ on $\R^3$ the isotropy of $b\neq0=a$ is $\g_b=\mathfrak{sl}_2\ltimes\R^2$, which is not Borel.
Alternatively, we can observe that the isotropy representation of $\g_b$ on $\R^3=\g/\g_b$ has no invariant 2-plane
(and only 1 invariant line). Thus this case is impossible.

As mentioned before, the case $\dim\g_a=7$ is impossible too, so consider the last possibility $\dim\g_a=6$. There is only one subalgebra of dimension 6 up to an outer automorphism 
(or two up to an inner automorphism): the maximal parabolic subalgebra $\mathfrak{p}_1$ 
(or $\mathfrak{p}_2$). Let $\mathcal{O}_a$ be the local 2-dimensional orbit through
$a$, it has homogeneous geometry of the type $SL_3/P_1$. For every point $b\in\mathcal{O}_a$
the isotropy in $\mathfrak{p}_1$ contains a simple Lie subalgebra $\mathfrak{h}_b\simeq\mathfrak{sl}_2$, 
analytically depending on $b$, which, by an application of the Guillemin-Sternberg theorem, is linearizable. 
Thus we get a curve $\ell_b$ of fixed points of $\mathfrak{h}_b$ passing through $b\in\mathcal{O}_a$. 
The union $\cup_{b\in\mathcal{O}_a}\ell_b$ is an open set
and hence contains a regular point $c\in\ell_b$ arbitrary close to $a$. The algebra $\mathfrak{h}_b$ belongs to the 
isotropy subalgebra $\g_c$ and so should leave invariant a 2-plane $\Pi_c\subset T_c\E$ split in two invariant lines. 
However there is only one invariant plane on which $\mathfrak{h}_b$ acts via the standard representation, 
so no invariant lines exist. 

This contradiction shows that the points in a neighborhood of $a$ are regular and hence 
all of them are flat with respect to the canonical Cartan connection; 
in other words, the ODE is locally trivializable, as required.

To prove the first part of the theorem it suffices to demonstrate that every integer between 0 and 6 is realizable
as the symmetry dimension of some ODE. Dimensions 0 to 3 are realizable via ODEs of principal type \cite{T}, so we
restrict to dimensions 6, 5, 4. Below are the examples realizing these dimensions (we indicate only the symmetry fields on $J^0(\R)=\R^2(x,y)$, their prolongations to $J^2(\R)$ are straightforward).
 \begin{gather*}
y\,y''=2(y')^2:\quad
\g=\langle \p_x,x\p_x,y\p_y,x^2\p_x-xy\p_y,y^2\p_y,xy^2\p_y\rangle,\\
xy\,y''=2y'(xy'+y):\quad
\g=\langle x\p_x,y\p_y,y^2\p_y,x^3y^2\p_y,x^4\p_x-3x^3y\p_y\rangle,\\
y\,y''=(y')^2:\quad
\g=\langle \p_x,x\p_x,y\p_y,xy\p_y\rangle.
 \end{gather*}
This finishes the proof. \qed

 \begin{rk}
The subalgebra of dimension 6 in $\mathfrak{sl}_3$ is unique up to outer automorphisms.
However the submaximal symmetric model is not unique. For every integer $k>1$ the ODE
$y\,y''=k(y')^2$ has 6-dimensional symmetry algebra 
$\g=\langle \p_x,x\p_x,y\p_y,(k-1)x^2\p_x-xy\p_y,y^k\p_y,xy^k\p_y\rangle$. 
The form of vector fields shows that these models are not mutually equivalent.
This is similar to the effect observed in \cite{IK2} for CR-structures.

Note also that the above ODE is of Fuchsian type $y''=-ky'/x$ after the 
interchange of variables $x\leftrightarrow y$. 
 \end{rk}


 \textsc{Acknowledgement.} The symbolic package \textsc{Maple} (DifferentialGeometry) was used
in computing symmetries.


\end{document}